\documentclass[11pt]{article}

\usepackage[dvips]{graphicx}
\usepackage{amssymb}
\usepackage{amsmath}

\usepackage{caption}
\usepackage{hyperref}
\usepackage{arcs}
\usepackage {xcolor} 
\usepackage{subfigure}
\usepackage{authblk}
\usepackage{calrsfs}

\newtheorem{theorem}{Theorem}

\newtheorem{Lemma}{Lemma}
\newtheorem{Remark}{Remark}

\usepackage{graphicx} % Required for inserting images

\title{Unisolvence of randomized MultiQuadric Kansa collocation  for convection-diffusion with \\mixed boundary conditions}

\author[1]{M. Mohammadi}
\affil[1]{Kharazmi University, Tehran}
\author[2]{A. Sommariva} 
\author[3]{M. Vianello\footnote{Corresponding author: marcov@math.unipd.it}}
\affil[2,3]{University of Padova, Italy}

\date{\today}

\begin{document}

\maketitle

\begin{abstract}
We make a further step in the open problem of unisolvence for unsymmetric
Kansa collocation, proving that the MultiQuadric Kansa method with fixed collocation points and random fictitious centers is almost surely unisolvent, for stationary convection-diffusion equations with mixed boundary conditions on general domains. 
For the purpose of illustration, the method is applied  in 2D with 
fictitious centers that are local random perturbations of predetermined collocation points.
\end{abstract}


% \end{frontmatter}

\section{Introduction}

In the recent paper \cite{MSV24} a further step has been made in the open problem of unisolvence for unsymmetric Kansa collocation, proving that collocation matrices for elliptic equations by Polyharmonic Splines (without polynomial addition) with fixed collocation points and random fictitious centers are almost surely invertible. Domains and boundary conditions are general (mixed type). The proving technique (by induction on determinants) is based on the fact that Polyharmonic Splines are real analytic but have a singularity at the center, that can be exploited to prove the key  property of linear independence of the functions involved in collocation.  

In the present note we extend such unisolvence result to the case of Kansa collocation by Multiquadrics, widely studied and applied after the pioneering work of E.J. Kansa \cite{K86,K90}; cf. e.g. \cite{C12,F07,LOS06,SW06} with the references therein. The proof is more difficult with respect to Polyharmonic Splines, since MultiQuadrics  are everywhere analytic so real singularities cannot be exploited to prove linear independence of the involved functions. Indeed, we have to resort to a complex embedding in order to to exploit the presence of complex singularities. 

We prove the result for stationary convection-diffusion equations on general domains with mixed Dirichlet-Neumann boundary conditions. A key aspect is that the centers are kept distinct from the collocation points and randomly chosen, 
so that the framework is different from \cite{CDRDASV24,DASV24} where the classical (but tricky) case of random collocation points coinciding with the centers was considered, for the Poisson equation with  purely Dirichlet boundary conditions. Indeed, the present framework allows to include in a simple way the Neumann conditions. 

Moreover, the fictitious centers are any continuous random vector, so that in practice they can be chosen as local random perturbations of both the interior and the boundary fixed collocation points, ensuring in any case almost sure unisolvence of the collocation process. The collocation points can be taken deterministically, for example with a uniform or quasi-uniform distribution in the interior and on the boundary of the domain. Determining sufficient conditions for unisolvence is relevant, since it is well-known after Hon and Schaback \cite{HS01} that it could not hold in the purely deterministic case, but there is still a substantial lack of theoretical results on this subject.

The paper is organized as follows. In Section 2 we give the main theoretical result, based on a quite general lemma on unisolvence of interpolation in analytic function spaces. In Section 3 we present a couple of numerical examples, showing the practical applicability of randomized Kansa collocation with fictitious centers, chosen as local random perturbations of predetermined collocation points.

\section{Unisolvence of MQ Kansa collocation}
We consider stationary convection-diffusion equations with  constants coefficients and mixed boundary conditions 
\begin{equation}\label{cd}
\left\{
\begin{array}{l}
\mathcal{L}u(P)=\Delta u(P)+\langle \nabla u(P),\vec{v}(P)\rangle=f(P)\;,\;P\in \Omega\subset \mathbb{R}^d\;,
\\ \\
\mathcal{B}u(P)=\chi_{\Gamma_1}(P)u(P)+\chi_{\Gamma_2}(P)\partial_\nu u(P)=g(P)\;,\;P\in \partial \Omega\;,
\end{array} 
\right .
\end{equation}
where $\Omega$ is a bounded domain (connected open set), 
$P=(x_1,\dots,x_d)$, $\Delta=\partial^2_{x_1}+\dots +\partial^2_{x_d}$ is the Laplacian, 
$\nabla=(\partial_{x_1},\dots,\partial_{x_d})$ denotes the gradient, $\vec{v}(P)$ a velocity field and $\langle \cdot,\cdot\rangle$ the inner product in $\mathbb{R}^d$,  
$\partial_\nu=\langle \nabla,\vec{\nu}\rangle$ is the normal derivative at a boundary point, and $\chi$ denotes the indicator function. The boundary is indeed splitted in two disjoint portions, namely $\partial \Omega=\Gamma_1\cup \Gamma_2$. If $\Gamma_2=\emptyset$
or $\Gamma_1=\emptyset$ we recover purely Dirichlet or purely Neumann conditions, respectively. 
Observe that for notational simplicity we have taken the diffusion coefficient equal to 1, with no loss of generality since otherwise we can absorb it in $\vec{v}$ and $f$. As known, equation (\ref{cd}) models the steady state of convection-diffusion with an incompressible flow. 
We do not make here any restrictive assumption on the domain $\Omega$ and on the functions $\vec{v}$, $f$ and $g$, except for those ensuring well-posedness and sufficient regularity of the solution (for example that the domain has a Lipschitz boundary, cf. e.g. \cite{S98} with the references therein).

We study the discretization of the convection-diffusion problem above by unsymmetric Kansa collocation using MQ (MultiQuadric) RBF (Radial Basis Functions) of the form $\{\phi_{C_j}(P)\}$, $1\leq j\leq N$, 
\begin{equation} \label{phiA}
\phi_C(P)=\phi(\|P-C\|)\;,\;\;\phi(r)=\sqrt{1+({\varepsilon}r)^2}\;,\;r\geq 0\;,
\end{equation}
where $C=(c_1,\dots,c_d)$ is the RBF center and $\|\cdot\|$ the Euclidean norm. As known, the so-called ``shape parameter''  $\varepsilon>0$ can be used to control
the trade-off between conditioning and accuracy; cf. e.g. \cite{C12,F07,LS23} with the references therein. 

The collocation points will be fixed, whereas the centers will be chosen as a random vector. The approach where centers are distinct from the collocation points is known as collocation by ``fictitious centers'' in the literature, differently from the classical method, originally  proposed in the pioneering work by E.J. Kansa \cite{K86,K90}, where centers and collocation points coincide. Methods based on fictitious centers are an active research subfield in  the literature on Kansa collocation (essentially in the least squares framework), 
cf. e.g. \cite{CKD19,CKA21,ZLHW22}. The possibility of taking separate collocation and center points allows more flexibility, both from the theoretical as well as the computational point of view. In particular, we will be able to prove almost sure unisolvence of the discretized problem. 

Seeking a solution of the form $u_N(P)=\sum_{j=1}^N{a_j\phi_{C_j}(P)}$ we get the linear system 
\begin{equation} \label{kansa-syst}
K_N\left(\begin{array} {c}
a_1\\ \vdots\\ a_N\\
\end{array} \right)=\left(\begin{array} {c}
f(P_i)\\ \\ g(P_k)\\
\end{array} \right)\;,\;\;K_N=\left[\begin{array} {c}
\mathcal{L} \phi_{C_j}(P_i)\\ \\
\mathcal{B} \phi_{C_j}(P_k)
\end{array} \right]\in \mathbb{R}^{N\times N}\;,
\end{equation}
where $1\leq i\leq N_I$, $N_I+1\leq k\leq N$, $1\leq j\leq N$, 
$C_1,\dots,C_N$ are the centers, $\{P_1,\dots,P_{N_I}\}\subset \Omega$ are $N_I$ distinct internal collocation points, and $\{P_{N_I+1},\dots,P_N\}\subset \partial\Omega$ are $N_B=N-N_I$ distinct boundary collocation points.

Observe that taking the partial derivatives with respect to the $P=(x_1,\dots,x_d)$ variables, we easily get (cf. e.g. \cite{F07})
\begin{equation} \label{lapl}
\nabla \phi_C(P)=
(P-C)\,\frac{\phi'(r)}{r}\;,\;\;\Delta\phi_C(P)=\phi''(r)+(d-1)\,\frac{\phi'(r)}{r}\;,\;\;r=\|P-C\|\;,
\end{equation}
so that 
\begin{equation} \label{L}
\mathcal{L} \phi_{C}(P)=\phi''(r)+\left(d-1+\langle P-C,\vec{v}(P)\rangle\right)\,\frac{\phi'(r)}{r}\;,\;\;P\in \Omega\;,
\end{equation}
\begin{equation} \label{B}
\mathcal{B} \phi_{C}(P)=\chi_{\Gamma_1}(P)\phi(r)+\chi_{\Gamma_2}(P)\langle P-C,\vec{\nu}(P)\rangle\,\frac{\phi'(r)}{r}\;,\;\;P\in \partial\Omega\;.
\end{equation}
In the particular case of MQ we have
\begin{equation} \label{derivMQ}
\frac{\phi'(r)}{r}=\varepsilon^2(1+(\varepsilon r)^2)^{-1/2}\;,\;\;\;\phi''(r)=-\varepsilon^4r^2(1+(\varepsilon r)^2)^{-3/2}+\varepsilon^2(1+(\varepsilon r)^2)^{-1/2}\;.
\end{equation}

We prove now a preliminary lemma on interpolation by analytic functions, that will be relevant below. 
\begin{Lemma}
Let $A\subseteq \mathbb{R}^d$, be an open connected set and $\{f_j\}_{1\leq j\leq N}$ be linearly independent real analytic functions in $A$. 

Then the set of non-unisolvent $N$-uples for interpolation in $\mbox{span}\{f_1,\dots,f_N\}$ has null Lebesgue measure in $A^N$. 
\end{Lemma}
\noindent
{\bf Proof.}
Consider determinant of the interpolation matrix   $$D(P_1,\dots,P_N)=\mbox{det}\left([f_j(P_i)]_{1\leq i,j\leq N}\right)\;,$$ 
as a function of $(P_1,\dots,P_N)\in A^N$. Such a function is analytic in $A^N$ since analytic functions form an algebra. Notice that $A^N$ is open and connected, being a product of open connected sets (cf. e.g. \cite{Mu14}). 
By a known general result on interpolation by linearly independent continuous functions  (cf. \cite{LOS06}), there exist $N$-uples in $A^N$ such that $D(P_1,\dots,P_N)$ does not vanish. 
Hence $D$ is not identically zero in $A^N$. In view of a fundamental theorem in the theory of real analytic functions (cf. e.g. \cite{M20}), then the zero set of $D$ in the open connected set $A^N$ has null Lebesgue measure. \hspace{0.2cm}  $\square$
\vskip0.2cm 

We are now ready to state and prove the following
\begin{theorem}
Let $K_N$ be the MQ Kansa collocation matrix in (\ref{kansa-syst}) for the convection-diffusion problem (\ref{cd}), where $\{P_i\,,\,1 \leq i\leq N_I\}\subset \Omega$ and $\{P_k\,,\,N_I+1\leq k\leq N\}\subset \partial\Omega$ are any two fixed sets of distinct internal and boundary collocation points, respectively, and $X=(C_1,\dots,C_N)$ a continuous random vector with probability density $\sigma(X)\in L^1(\mathbb{R}^{dN})$.

Then the matrix $K_N$ is almost surely nonsingular.
\end{theorem}
\vskip0.2cm 
\noindent 
{\bf Proof.} The key observation is that, once fixed the set of distinct collocation points $\{P_i\}$, the matrix $K_N$ can be seen as the interpolation matrix at the points $C_1,\dots,C_N$, with the functions 
$$
f_1(C)=\mathcal{L} \phi_{C}(P_1),\dots,f_{N_I}(C)=\mathcal{L} \phi_{C}(P_{N_I}),
$$
\begin{equation} \label{f1N}
f_{N_I+1}(C)=\mathcal{B} \phi_{C}(P_{N_I+1}),\dots,f_{N}(C)=\mathcal{B} \phi_{C}(P_{N})\;.
\end{equation}
Such functions are real analytic in $A=\mathbb{R}^{dN}$, in view of (\ref{lapl})-(\ref{derivMQ}) and the analyticity in $r\in \mathbb{R}$  of the univariate function $(1+(\varepsilon r)^2)^s$, $s\in \mathbb{R}$.
In order to apply Lemma 2.1, we have to prove that the functions $f_1(C),\dots,f_N(C)$ are linear independent. 

Now, assume that they were dependent. Then, there is an everywhere vanishing  linear combination, $F(C)=\sum_{j=1}^N{\alpha_jf_j(C)}\equiv 0$, with $\alpha_{\ell}\neq 0$ for some $\ell$. Take the line $C(t)=P_{\ell}+tU$ with any fixed unit vector $U=(u_1,\dots,u_d)$, then the univariate analytic function $F(C(t))$ is identically zero in $\mathbb{R}$ and thus its complex extension $F(C(z))$ is identically zero in $\mathbb{C}$. Notice that the functions 
$(1+\varepsilon^2\|P_j-C(z)\|)^s=(1+\varepsilon^2\|P_j-P_\ell-zU\|^2)^s$ for $s=1/2,-1/2,-3/2$, appearing in the complex extension of the functions $f_j(C(z))$, correspond to the branch of the fractional powers which is positive on the real positive axis. Moreover, 
$\|P_j-P_\ell-zU\|^2$ has to be seen as the complex extension of the corresponding real function, hence not the complex 2-norm but the sum of the squares of the complex components. 

Then $(1+\varepsilon^2\|P_\ell-C(z)\|^2)^s=(1+\varepsilon^2z^2)^s$ presents two branching points at $z=\pm i/\varepsilon$, whereas for $j\neq \ell$ the functions 
$(1+\varepsilon^2\|P_j-P_\ell-zU\|^2)^s$ are analytic at $z=\pm i/\varepsilon$, since the complex numbers $$
1+\varepsilon^2\|P_j-P_\ell-(\pm i/\varepsilon)U\|^2=1+\varepsilon^2\sum_{h=1}^d(P_j-P_\ell\mp iU/\varepsilon)_h^2
$$
$$
=1+\varepsilon^2\sum_{h=1}^d[(P_j-P_\ell)_h^2\mp 2i(P_j-P_\ell)_hu_h/\varepsilon-u_h^2/\varepsilon^2] 
$$
$$
=\varepsilon^2\sum_{h=1}^d(P_j-P_\ell)_h^2\mp 2i\varepsilon\sum_{h=1}^d(P_j-P_\ell)_hu_h
$$ 
have positive real part. But $F(C(z))\equiv 0$ means that $f_\ell(C(z))$ is a linear combination of the functions 
$f_j(C(z))$, $j\neq \ell$, and this gives immediately a contradiction, since the latter are analytic at the branching points present in 
$f_\ell(C(z))$ by (\ref{L})-(\ref{derivMQ}). 

At this point we can apply Lemma 2.1, obtaining that the set of $N$-uples of centers  $X=(C_1,\dots,C_N)$ for which the collocation matrix $K_N$ is singular, has null Lebesgue measure in $\mathbb{R}^{dN}$. Consequently, it has null measure with respect to 
any absolutely continuous measure with respect to the Lebesgue measure, and hence the matrix $K_N$ is almost surely nonsingular for any distribution of centers by 
a continuous probability measure with density $\sigma(X)\in L^1(\mathbb{R}^{dN})$. \hspace{0.2cm} $\square$

\begin{Remark}
It is worth observing that the possible center distribution is more general than that assumed in \cite{MSV24}, where the fictitious centers 
are a sequence of i.i.d. (indipendent identically distributed) random points. 
Indeed in the present framework the multivariate probability density 
may not even be a product density. 
\end{Remark}

\section{Numerical examples} 
Though the main purpose of the present work is theoretical, making a further step within the theoretical open problem of Kansa collocation unisolvence, the method of random fictitious centers can be conveniently adopted with suitable cautions and tricks. 
Indeed, while the random center distribution in Theorem 2.1 is quite general and the random centers could be placed in principle anywhere, still ensuring almost sure unisolvence, in practice the method works much better with centers located near the collocation points. On the other hand, it is well-known that the MQ collocation matrices can be severely ill-conditioned and a specialized literature exists on different approaches to cope with ill-conditioning, such as for example shape parameter optimization, extended precision arithmetic, RBF-QR method; cf., e.g., \cite{F07,KH23,LLHF13,LS23} with the references therein.

For the mere purpose of illustration we present some simple numerical examples, concerning the solution of convection-diffusion equations with mixed boundary conditions on a square, by randomized Kansa collocation with MultiQuadrics, implemented in Matlab. 
We consider the
convection-diffusion problem (\ref{cd}) on $\Omega=(0,1)^2$.
In all the test problems, the set of fixed collocation points $\mathcal{C}=\{P_j\}_{1\leq j\leq N}$ is a uniform grid on the square, in lexicographic order. Using a Matlab notation, the random fictitious centers $$X=\mathcal{C}+(2* \mbox{rand}(N,1)-1)*\delta$$ are obtained by local random perturbation of the collocation points via additive uniformly distributed random points in $(-\delta,\delta)^2$; see Fig. 1. 
The accuracy is measured by the geometric mean of Root Mean Square Errors 
$\mbox{RMSE}_{av}=\exp\left(\frac{1}{m} \,\sum_{l=1}^m{\log_{10} \left(\sqrt{\sum_j{(u_j-\tilde{u}_{j,l})^2/N}}\right)}\right)$
obtained by $m$ random centers arrays $\{X_l\}$, $l=1,\ldots,m$, where $u_j$ and $\tilde{u}_{j,l}$ are the exact and approximate
solutions at the collocation node {$P_j$}, respectively. In the tests, we have run $m=100$ trials.

We have imposed mixed-type boundary conditions in (\ref{cd}), by the splitting 
$$\Gamma_1=\{x_1=0,~ 0\leq x_2\leq 1\}\cup  \{x_1=1,~ 0\leq x_2\leq 1\}\;,
$$
$$
\Gamma_2=\{x_2=0,~ 0< x_1< 1\}\cup  \{x_2=1,~ 0< x_1< 1\}\;.
$$
The right-hand sides $f$ and $g$ are defined by selecting as reference solution 
$u(x_1,x_2)=\sin(2\pi x_1)+\cos(2\pi x_2)$. 

The numerical results are collected in Table 1, where we have taken different values of $\delta$ and different convection velocities, which make the problem ranging from pure diffusion to mildly convection-dominated instances, and we have found heuristically a value of the shape parameter $\varepsilon$, which roughly minimizes the errors. For strongly convection-dominated problems with high P\'{e}clet number, more specific discretization techniques should be adopted, that go beyond the scope of the present paper; cf. e.g. \cite{WH24} with the references therein. 

We can observe that for the smallest values of $\delta$ the errors approach the size of those corresponding to classical collocation with  $X=\mathcal{C}$, which can be considered a limit case, whose unisolvence is however not covered by the present theory. For larger values 
of $\delta$, e.g. $\delta=0.1$, the errors are not satisfactory, which could be due to the fact that less centers fall near the boundary.

\begin{figure}[!]
%  % Requires \usepackage{graphicx}
 \center{ \includegraphics[scale=0.40]{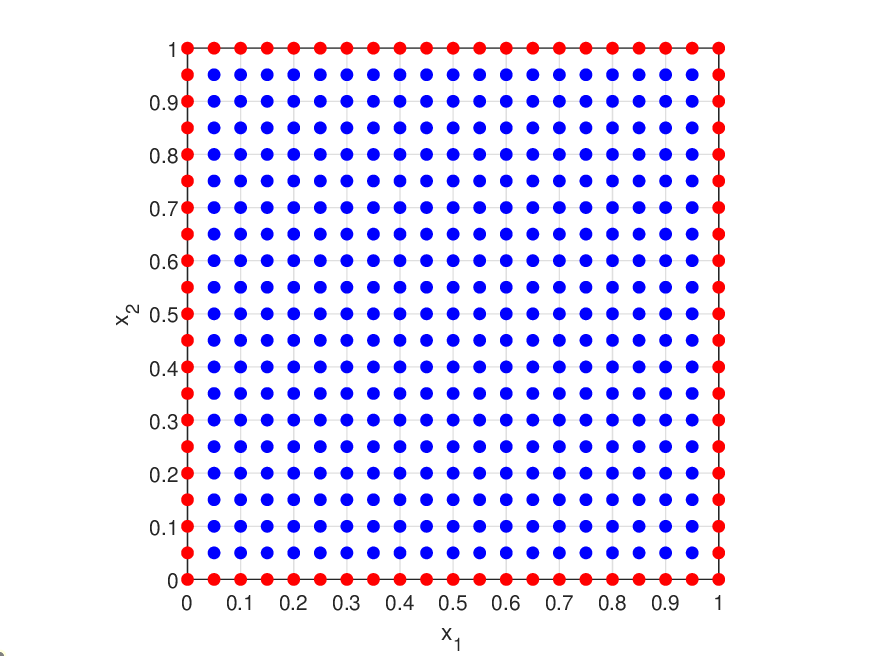} \includegraphics[scale=0.40 ]{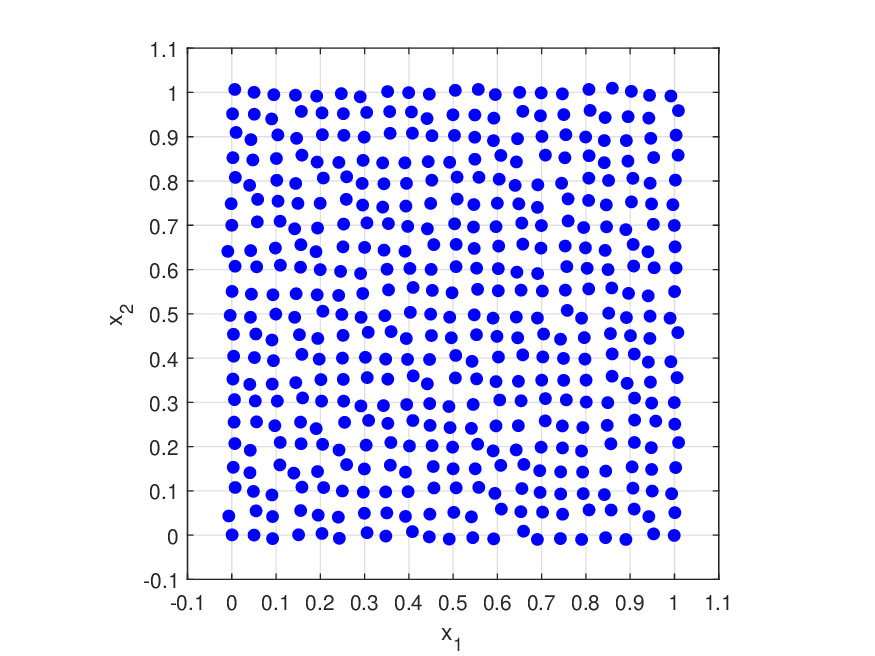}}
    \caption {$441$ collocation grid points (left) and the random fictitious centers distribution (right) for $\delta=0.01$.}\label{eqpointstest1}
    \end{figure}

\begin{table}[!]
  \centering  
  \caption{RMSE geometric mean over 100 trials of $N$ random fictitious centers, with different convection velocities $\vec{v}$ and perturbation radius $\delta$ (shape parameter $\varepsilon=2.5$).} 
  \label{t1}  
  \vspace{5pt} % Add vertical space before the table  
  \scriptsize{  
    \begin{tabular}{|c|c|c|c|c||c|c|c|c|}  
      \hline 
      & \multicolumn{4}{c||}{$\vec{v}=(0,0)$} & \multicolumn{4}{c|}{$\vec{v}=(1,1)$} \\
      \hline  
      $N$ & $\;\;\delta=0.1\;\;$ &  $\;\delta=0.01\;$ & $\delta=0.001$ & $\;\;\;\delta=0\;\;\;$ & $\;\;\delta=0.1\;\;$ &  $\;\delta=0.01\;$ & $\delta=0.001$ & $\;\;\;\delta=0\;\;\;$\\
      \hline  
      $121$ &  2.9e-01 &   6.9e-02 &   6.8e-02  &  6.9e-02 & 3.4e-01 &   7.1e-02  &  7.3e-02  & 7.2e-02 \\
      \hline
      $441$  &  4.7e-02  &  7.9e-03  &  1.8e-03  &  1.4e-03 &  6.2e-02  &  8.2e-03  &  1.7e-03 &   1.5e-03 \\
      \hline
      $961$  & 1.4e-03  &  9.6e-04 &   4.5e-04 &   3.4e-05 & 1.7e-03  &  9.1e-04  &  4.5e-04  &  3.5e-05 \\
      \hline
      $1681$  & 4.3e-05  &  9.8e-06 &   1.0e-05  &  7.5e-06 & 4.7e-05  &  1.2e-05  &  1.2e-05 &   6.3e-06 \\
      \hline  
      
      \hline 
      & \multicolumn{4}{c||}{$\vec{v}=(1,100)$} & \multicolumn{4}{c|}{$\vec{v}=(100,100)$} \\
      \hline  
      $N$ & $\;\;\delta=0.1\;\;$ &  $\;\delta=0.01\;$ & $\delta=0.001$ & $\;\;\;\delta=0\;\;\;$ & $\;\;\delta=0.1\;\;$ &  $\;\delta=0.01\;$ & $\delta=0.001$ & $\;\;\;\delta=0\;\;\;$\\
      \hline  
      $121$ & 6.1e+00  &  5.1e+00  &  3.1e+00 &   3.1e+00 &  1.6e-01  &  8.3e-02 &   4.7e-01  &  3.9e-01\\
      \hline
      $441$  & 2.4e+00  &  6.5e-01  &  1.1e-01 &   6.5e-02 &  1.4e-02   &  5.9e-03  &  3.7e-03 &   4.2e-03\\
      \hline
      $961$  & 7.7e-02   & 3.7e-02 &   2.4e-02  &  1.5e-03 &  8.5e-04   & 4.9e-04&    2.5e-04  &  1.2e-04 \\ 
      \hline
      $1681$  & 3.1e-03  &  8.1e-04  &  7.9e-04 &   2.1e-04  & 4.8e-05  &  2.3e-05  &  1.8e-05   & 1.9e-05 \\
      \hline  
      \end{tabular} 
  }  
  \vspace{5pt} % Add vertical space after the table  
\end{table}

\newpage
\noindent
{\footnotesize 
{\bf Acknowledgements.} 

Work partially supported by the DOR funds of the University of Padova, by the INdAM-GNCS 2024 Project “Kernel and polynomial methods for approximation and integration: theory and application software'' and 2025 Project Polynomials, Splines and Kernel Functions: from Numerical Approximation to Open-Source Software'' (A. Sommariva and M. Vianello), and by the INdAM-GNCS 2024 Project “Meshless Techniques and integro-differential equations: analysis and their application'' (M. Mohammadi). This research has been accomplished within the RITA ``Research ITalian network on Approximation", the SIMAI Activity Group ANA\&A, and the UMI Group TAA ``Approximation Theory and Applications".
\/}


\begin{thebibliography}{99}

\bibitem{CDRDASV24} R. Cavoretto, A. De Rossi, F. Dell'Accio, A. Sommariva, M. Vianello, Nonsingularity of unsymmetric Kansa matrices: random collocation by MultiQuadrics and Inverse MultiQuadrics, arXiv: 2403.18017.

\bibitem{CKD19} C. Chen, A. Karageorghis, F. Dou, A novel RBF collocation method using fictitious centres, Appl. Math. Lett. 101 (2019), 106069.

\bibitem{CKA21} C.S. Chen, A. Karageorghis, L. Amuzu,  
Kansa RBF collocation method with auxiliary boundary
centres for high order BVPs, J. Comput. Appl. Math. 398 (2021), 113680.

\bibitem{C12} A.H.-D Cheng, Multiquadric and its shape parameter—A numerical investigation
of error estimate, condition number, and round-off error by arbitrary
precision computation, Eng. Anal. Bound. Elem. 36 (2012), 220--239.

\bibitem{DASV24} F. Dell'Accio, A. Sommariva, M. Vianello, Unisolvence of random Kansa collocation by Thin-Plate Splines for the Poisson equation, 
Eng. Anal. Bound. Elem. 165 (2024), 105773.

\bibitem{F07} G.E. Fasshauer, Meshfree Approximation Methods with Matlab, Interdisciplinary Mathematical
Sciences, Vol. 6, World Scientific, 2007.

\bibitem{HS01}
Y.C. Hon, R. Schaback, On unsymmetric collocation by radial basis functions, Appl. Math. Comput. 119 (2001), 177--186.

\bibitem{K86} E.J. Kansa,  Application of Hardy’s multiquadric interpolation to hydrodynamics, Lawrence Livermore National Lab. preprint, CA (USA), 1985 (Proc. 1986 Simul. Conf., San Diego, pp. 111–117).

\bibitem{K90} E.J. Kansa, Multiquadrics - A scattered data approximation scheme with applications to computational fluid-dynamics - II solutions to parabolic, hyperbolic and elliptic partial differential equations, Comput. Math. Appl. 19 (1990), 147--161.

\bibitem{KH23} E.J. Kansa, P. Holoborodko, On the ill-conditioned nature of $C^ \infty$ RBF strong collocation, 
Eng. Anal. Bound. Elem. 78 (2017), 26--30.

\bibitem{LLHF13} E. Larsson, E. Lehto, A. Heryudono, B. Fornberg, Stable computation of differentiation matrices and scattered node stencils based on Gaussian radial basis functions, SIAM J. Sci. Comput. 35 (2013), A2096--A2119.

\bibitem{LS23} E. Larsson, R. Schaback, Scaling of radial basis functions, IMA Journal of
Numerical Analysis (2023), drad035.

\bibitem{LOS06} L. Ling, R. Opfer, R. Schaback, Results on meshless collocation techniques, Eng. Anal. Bound. Elem. 30 (2006), 247--253.

\bibitem{M20} B.S. Mityagin, The Zero Set of a Real Analytic Function, Math. Notes 107 (2020), 529--530.

\bibitem{MSV24} M. Mohammadi, A. Sommariva, M. Vianello, Unisolvence of Kansa collocation for elliptic equations by polyharmonic splines with random fictitious centers, arXiv: 2410.03279. 

\bibitem{Mu14} J. Munkres, Topology, 2nd edition, Pearson, 2014.

\bibitem{S98} S. Savar\`{e}, Regularity Results for Elliptic Equations in
Lipschitz Domains, J. Funct. Anal. 152 (1998), 176--201.

\bibitem{SW06} R. Schaback, H. Wendland, 
Kernel techniques: From machine learning to meshless methods, Acta Numer. 15 (2006,  543--639.

\bibitem{WH24} J. Wang, M. Hillman, Upwind reproducing kernel collocation method for convection-dominated problems, Comput. Methods Appl. Mech. Eng. 420 (2024), 116711. 

\bibitem{ZLHW22} H. Zheng, X. Lai, A. Hong, X. Wei, 
A Novel RBF Collocation Method Using Fictitious Centre Nodes for Elasticity Problems, MDPI Mathematics 10(19) (2022), 3711.

\end{thebibliography}
\end{document}